# HIGH-DIMENSIONAL VARIABLE SELECTION

By Larry Wasserman and Kathryn Roeder[1]

*Carnegie Mellon University*

This paper explores the following question: what kind of statistical guarantees can be given when doing variable selection in high-dimensional models? In particular, we look at the error rates and power of some multi-stage regression methods. In the first stage we fit a set of candidate models. In the second stage we select one model by cross-validation. In the third stage we use hypothesis testing to eliminate some variables. We refer to the first two stages as "screening" and the last stage as "cleaning." We consider three screening methods: the lasso, marginal regression, and forward stepwise regression. Our method gives consistent variable selection under certain conditions.

**1. Introduction.** Several methods have been developed lately for high-dimensional linear regression such as the lasso [Tibshirani (1996)], Lars [Efron et al. (2004)] and boosting [Bühlmann (2006)]. There are at least two different goals when using these methods. The first is to find models with good prediction error. The second is to estimate the true "sparsity pattern," that is, the set of covariates with nonzero regression coefficients. These goals are quite different and this paper will deal with the second goal. (Some discussion of prediction is in the Appendix.) Other papers on this topic include Meinshausen and Bühlmann (2006), Candes and Tao (2007), Wainwright (2006), Zhao and Yu (2006), Zou (2006), Fan and Lv (2008), Meinshausen and Yu (2008), Tropp (2004, 2006), Donoho (2006) and Zhang and Huang (2006). In particular, the current paper builds on ideas in Meinshausen and Yu (2008) and Meinshausen (2007).

Let $(X_1, Y_1), \ldots, (X_n, Y_n)$ be i.i.d. observations from the regression model

$$Y_i = X_i^T \beta + \varepsilon_i, \tag{1}$$

Received June 2008; revised August 2008.
[1]Supported by NIH Grant MH057881.

where $\varepsilon \sim N(0, \sigma^2)$, $X_i = (X_{i1}, \ldots, X_{ip})^T \in \mathbb{R}^p$ and $p = p_n > n$. Let $X$ be the $n \times p$ design matrix with $j$th column $X_{\bullet j} = (X_{1j}, \ldots, X_{nj})^T$ and let $Y = (Y_1, \ldots, Y_n)^T$. Let

$$D = \{j : \beta_j \neq 0\}$$

be the set of covariates with nonzero regression coefficients. Without loss of generality, assume that $D = \{1, \ldots, s\}$ for some $s$. A variable selection procedure $\widehat{D}_n$ maps the data into subsets of $S = \{1, \ldots, p\}$.

The main goal of this paper is to derive a procedure $\widehat{D}_n$ such that

(2) $$\limsup_{n \to \infty} \mathbb{P}(\widehat{D}_n \subset D) \geq 1 - \alpha,$$

that is, the asymptotic type I error is no more than $\alpha$. Note that throughout the paper we use $\subset$ to denote nonstrict set-inclusion. Moreover, we want $\widehat{D}_n$ to have nontrivial power. Meinshausen and Bühlmann (2006) control a different error measure. Their method guarantees $\limsup_{n \to \infty} \mathbb{P}(\widehat{D}_n \cap V \neq \varnothing) \leq \alpha$ where $V$ is the set of variables not connected to $Y$ by any path in an undirected graph.

Our procedure involves three stages. In stage I we fit a suite of candidate models, each model depending on a tuning parameter $\lambda$,

$$\mathcal{S} = \{\widehat{S}_n(\lambda) : \lambda \in \Lambda\}.$$

In stage II we select one of those models $\widehat{S}_n$ using cross-validation to select $\widehat{\lambda}$. In stage III we eliminate some variables by hypothesis testing. Schematically,

$$\underbrace{\text{data} \xrightarrow{\text{stage I}} \mathcal{S} \xrightarrow{\text{stage II}} \widehat{S}_n}_{\text{screen}} \underbrace{\xrightarrow{\text{stage III}} \widehat{D}_n}_{\text{clean}}.$$

Genetic epidemiology provides a natural setting for applying screen and clean. Typically, the number of subjects, $n$, is in the thousands, while $p$ ranges from tens of thousands to hundereds of thousands of genetic features. The number of genes exhibiting a detectable association with a trait is extremely small. Indeed, for type I diabetes only ten genes have exhibited a reproducible signal [Wellcome Trust (2007)]. Hence, it is natural to assume that the true model is sparse. A common experimental design involves a 2-stage sampling of data, with stages 1 and 2 corresponding to the screening and cleaning processes, respectively.

In stage 1 of a genetic association study, $n_1$ subjects are sampled and one or more traits such as bone mineral density are recorded. Each subject is also measured at $p$ locations on the chromosomes. These genetic covariates usually have two forms in the population due to variability at a single nucleotide and hence are called single nucleotide polymorphisms (SNPs). The distinct forms are called alleles. Each covariate takes on a value (0, 1 or 2) indicating the number of copies of the less common allele observed. For a well-designed genetic study, individual SNPs are nearly uncorrelated unless they are physically located in very close proximity. This feature makes it much easier to draw causal inferences about the relationship between SNPs and quantitative traits. It is standard in the field to infer that an association discovered between a SNP and a quantitative trait implies a causal genetic variant is physically located near the one exhibiting association. In stage 2, $n_2$ subjects are sampled at a subset of the SNPs assessed in stage 1. SNPs measured in stage 2 are often those that achieved a test statistic that exceeded a predetermined threshold of significance in stage 1. In essence, the two stage design pairs naturally with a screen and clean procedure.

For the screen and clean procedure, it is essential that $\widehat{S}_n$ has two properties as $n \to \infty$ as follows:

$$\text{(3)} \qquad \mathbb{P}(D \subset \widehat{S}_n) \to 1$$

and

$$\text{(4)} \qquad |\widehat{S}_n| = o_P(n),$$

where $|M|$ denotes the number of elements in a set $M$. Condition (3) ensures the validity of the test in stage III while condition (4) ensures that the power of the test is not too small. Without condition (3), the hypothesis test in stage III would be biased. We will see that the power goes to 1, so taking $\alpha = \alpha_n \to 0$ implies consistency: $\mathbb{P}(\widehat{D}_n = D) \to 1$. For fixed $\alpha$, the method also produces a confidence sandwich for $D$, namely,

$$\liminf_{n \to \infty} \mathbb{P}(\widehat{D}_n \subset D \subset \widehat{S}_n) \geq 1 - \alpha.$$

To fit the suite of candidate models, we consider three methods. In method 1,

$$\widehat{S}_n(\lambda) = \{j : \widetilde{\beta}_j(\lambda) \neq 0\},$$

where $\widetilde{\beta}_j(\lambda)$ is the lasso estimator, the value of $\beta$ that minimizes

$$\sum_{i=1}^{n}(Y_i - X_i^T \beta)^2 + \lambda \sum_{j=1}^{p} |\beta_j|.$$

In method 2, take $\widehat{S}_n(\lambda)$ to be the set of variables chosen by forward stepwise regression after $\lambda$ steps. In method 3, marginal regression, we take

$$\widehat{S}_n = \{j : |\widehat{\mu}_j| > \lambda\},$$

where $\widehat{\mu}_j$ is the marginal regression coefficient from regressing $Y$ on $X_j$. (This is equivalent to ordering by the absolute $t$-statistics since we will assume that the covariates are standardized.) These three methods are very similar to basis pursuit, orthogonal matching pursuit and thresholding [see, e.g., Tropp (2004, 2006) and Donoho (2006)].

*Notation.* Let $\psi = \min_{j \in D} |\beta_j|$. Define the loss of any estimator $\widehat{\beta}$ by

$$L(\widehat{\beta}) = \frac{1}{n}(\widehat{\beta} - \beta)^T X^T X (\widehat{\beta} - \beta) = (\widehat{\beta} - \beta)^T \widehat{\Sigma}_n (\widehat{\beta} - \beta), \tag{5}$$

where $\widehat{\Sigma}_n = n^{-1} X^T X$. For convenience, when $\widehat{\beta} \equiv \widehat{\beta}(\lambda)$ depends on $\lambda$ we write $L(\lambda)$ instead of $L(\widehat{\beta}(\lambda))$. If $M \subset S$, let $X_M$ be the design matrix with columns $(X_{\bullet j} : j \in M)$ and let $\widehat{\beta}_M = (X_M^T X_M)^{-1} X_M^T Y$ denote the least-squares estimator, assuming it is well defined. Note that our use of $X_{\bullet j}$ differs from standard ANOVA notation. Write $X_\lambda$ instead of $X_M$ when $M = \widehat{S}_n(\lambda)$. When convenient, we extend $\widehat{\beta}_M$ to length $p$ by setting $\widehat{\beta}_M(j) = 0$, for $j \notin M$. We use the norms

$$\|v\| = \sqrt{\sum_j v_j^2}, \qquad \|v\|_1 = \sum_j |v_j| \quad \text{and} \quad \|v\|_\infty = \max_j |v_j|.$$

If $C$ is any square matrix, let $\phi(C)$ and $\Phi(C)$ denote the smallest and largest eigenvalues of $C$. Also, if $k$ is an integer, define

$$\phi_n(k) = \min_{M : |M| = k} \phi\left(\frac{1}{n} X_M^T X_M\right) \quad \text{and} \quad \Phi_n(k) = \max_{M : |M| = k} \Phi\left(\frac{1}{n} X_M^T X_M\right).$$

We will write $z_u$ for the upper quantile of a standard normal, so that $\mathbb{P}(Z > z_u) = u$ where $Z \sim N(0,1)$.

Our method will involve splitting the data randomly into three groups $\mathcal{D}_1$, $\mathcal{D}_2$ and $\mathcal{D}_3$. For ease of notation, assume the total sample size is $3n$ and that the sample size of each group is $n$.

*Summary of assumptions.* We will use the following assumptions throughout except in Section 8:

(A1) $Y_i = X_i^T \beta + \varepsilon_i$ where $\varepsilon_i \sim N(0, \sigma^2)$, for $i = 1, \ldots, n$.

(A2) The dimension $p_n$ of $X$ satisfies $p_n \to \infty$ and $p_n \leq c_1 e^{n^{c_2}}$, for some $c_1 > 0$ and $0 \leq c_2 < 1$.

(A3) $s \equiv |\{j : \beta_j \neq 0\}| = O(1)$ and $\psi = \min\{|\beta_j| : \beta_j \neq 0\} > 0$.

(A4) There exist positive constants $C_0, C_1$ and $\kappa$ such that $\mathbb{P}(\limsup_{n \to \infty} \Phi_n(n) \leq C_0) = 1$ and $\mathbb{P}(\liminf_{n \to \infty} \phi_n(C_1 \log n) \geq \kappa) = 1$. Also, $\mathbb{P}(\phi_n(n) > 0) = 1$, for all $n$.

(A5) The covariates are standardized: $\mathbb{E}(X_{ij}) = 0$ and $\mathbb{E}(X_{ij}^2) = 1$. Also, there exists $0 < B < \infty$ such that $\mathbb{P}(|X_{jk}| \leq B) = 1$.

For simplicity, we include no intercepts in the regressions. The assumptions can be weakened at the expense of more complicated proofs. In particular, we can let $s$ increase with $n$ and $\psi$ decrease with $n$. Similarly, the normality and constant variance assumptions can be relaxed.

**2. Error control.** Define the type I error rate $q(\widehat{D}_n) = \mathbb{P}(\widehat{D}_n \cap D^c \neq \varnothing)$ and the asymptotic error rate $\limsup_{n \to \infty} q(\widehat{D}_n)$. We define the power $\pi(\widehat{D}_n) = \mathbb{P}(D \subset \widehat{D}_n)$ and the average power

$$\pi_{\mathrm{av}} = \frac{1}{s} \sum_{j \in D} \mathbb{P}(j \in \widehat{D}_n).$$

It is well known that controlling the error rate is difficult for at least three reasons: correlation of covariates, high-dimensionality of the covariate and unfaithfulness (cancellations of correlations due to confounding). Let us briefly review these issues.

It is easy to construct examples where $q(\widehat{D}_n) \leq \alpha$ implies that $\pi(\widehat{D}_n) \approx \alpha$. Consider the following two models for random variables $Z = (Y, X_1, X_2)$:

| Model 1 | Model 2 |
|---|---|
| $X_1 \sim N(0,1),$ | $X_2 \sim N(0,1),$ |
| $Y = \psi X_1 + N(0,1),$ | $Y = \psi X_2 + N(0,1),$ |
| $X_2 = \rho X_1 + N(0,\tau^2).$ | $X_1 = \rho X_2 + N(0,\tau^2).$ |

Under models 1 and 2, the marginal distribution of $Z$ is $P_1 = N(\mathbf{0}, \Sigma_1)$ and $P_2 = N(\mathbf{0}, \Sigma_2)$, where

$$\Sigma_1 = \begin{pmatrix} \psi^2 + 1 & \psi & \rho\psi \\ \psi & 1 & \rho \\ \rho\psi & \rho & \rho^2 + \tau^2 \end{pmatrix}, \qquad \Sigma_2 = \begin{pmatrix} \psi^2 + 1 & \rho\psi & \psi \\ \rho\psi & \rho^2 + \tau^2 & \rho \\ \psi & \rho & 1 \end{pmatrix}.$$

Given any $\varepsilon > 0$, we can choose $\rho$ sufficiently close to 1 and $\tau$ sufficiently close to 0 such that $\Sigma_1$ and $\Sigma_2$ are as close as we like, and hence, $d(P_1^n, P_2^n) < \varepsilon$ where $d$ is total variation distance. It follows that

$$\mathbb{P}_2(2 \notin \widehat{D}) \geq \mathbb{P}_1(2 \notin \widehat{D}) - \varepsilon \geq 1 - \alpha - \varepsilon.$$

Thus, if $q \leq \alpha$, then the power is less than $\alpha + \varepsilon$.

Dimensionality is less of an issue thanks to recent methods. Most methods, including those in this paper, allow $p_n$ to grow exponentially. But all the methods require some restrictions on the number $s$ of nonzero $\beta_j$'s. In other words, some sparsity assumption is required. In this paper we take $s$ fixed and allow $p_n$ to grow.

False negatives can occur during screening due to cancellations of correlations. For example, the correlation between $Y$ and $X_1$ can be 0 even when $\beta_1$ is huge. This problem is called unfaithfulness in the causality literature [see Spirtes, Glymour and Scheines (2001) and Robins et al. (2003)]. False negatives during screening can lead to false positives during the second stage.

Let $\widehat{\mu}_j$ denote the regression coefficient from regressing $Y$ on $X_j$. Fix $j \leq s$ and note that

$$\mu_j \equiv \mathbb{E}(\widehat{\mu}_j) = \beta_j + \sum_{\substack{k \neq j \\ 1 \leq k \leq s}} \beta_k \rho_{kj},$$

where $\rho_{kj} = \operatorname{corr}(X_k, X_j)$. If

$$\sum_{\substack{k \neq j \\ 1 \leq k \leq s}} \beta_k \rho_{kj} \approx -\beta_j,$$

then $\mu_j \approx 0$ no matter how large $\beta_j$ is. This problem can occur even when $n$ is large and $p$ is small.

For example, suppose that $\beta = (10, -10, 0, 0)$ and that $\rho(X_i, X_j) = 0$ except that $\rho(X_1, X_2) = \rho(X_1, X_3) = \rho(X_2, X_4) = 1 - \varepsilon$, where $\varepsilon > 0$ is small. Then,

$$\beta = (10, -10, 0, 0), \quad \text{but } \mu \approx (0, 0, 10, -10).$$

Marginal regression is extremely susceptible to unfaithfulness. The lasso and forward stepwise, less so. However, unobserved covariates can induce unfaithfulness in all the methods.

**3. Loss and cross-validation.** Let $X_\lambda = (X_{\bullet j} : j \in \widehat{S}_n(\lambda))$ denote the design matrix corresponding to the covariates in $\widehat{S}_n(\lambda)$ and let $\widehat{\beta}(\lambda)$ be the least-squares estimator for the regression restricted to $\widehat{S}_n(\lambda)$, assuming the estimator is well defined. Hence, $\widehat{\beta}(\lambda) = (X_\lambda^T X_\lambda)^{-1} X_\lambda^T Y$. More generally, $\widehat{\beta}_M$ is the least-squares estimator for any subset of variables $M$. When convenient, we extend $\widehat{\beta}(\lambda)$ to length $p$ by setting $\widehat{\beta}_j(\lambda) = 0$, for $j \notin \widehat{S}_n(\lambda)$.

3.1. *Loss.* Now we record some properties of the loss function. The first part of the following lemma is essentially Lemma 3 of Meinshausen and Yu (2008).

LEMMA 3.1. *Let $\mathcal{M}_m^+ = \{M \subset S : |M| \leq m, D \subset M\}$. Then,*

(6) $$\mathbb{P}\left(\sup_{M \in \mathcal{M}_m^+} L(\widehat{\beta}_M) \leq \frac{4m \log p}{n \phi_n(m)}\right) \to 1.$$

*Let $\mathcal{M}_m^- = \{M \subset S : |M| \leq m, D \not\subset M\}$. Then,*

(7) $$\mathbb{P}\left(\inf_{M \in \mathcal{M}_m^-} L(\widehat{\beta}_M) \geq \psi^2 \phi_n(m+s)\right) \to 1.$$

3.2. *Cross-validation.* Recall that the data have been split into groups $\mathcal{D}_1$, $\mathcal{D}_2$ and $\mathcal{D}_3$ each of size $n$. Construct $\widehat{\beta}(\lambda)$ from $\mathcal{D}_1$ and let

$$\widehat{L}(\lambda) = \frac{1}{n} \sum_{X_i \in \mathcal{D}_2} (Y_i - X_i^T \widehat{\beta}(\lambda))^2. \tag{8}$$

We would like $\widehat{L}(\lambda)$ to order the models the same way as the true loss $L(\lambda)$ [defined after (5)]. This requires that, asymptotically, $\widehat{L}(\lambda) - L(\lambda) \approx \delta_n$, where $\delta_n$ does not involve $\lambda$. The following bounds will be useful. Note that $L(\lambda)$ and $\widehat{L}(\lambda)$ are both step functions that only change value when a variable enters or leaves the model.

THEOREM 3.2. *Suppose that $\max_{\lambda \in \Lambda_n} |\widehat{S}_n(\lambda)| \leq k_n$. Then, there exists a sequence of random variables $\delta_n = O_P(1)$ that do not depend on $\lambda$ or $X$, such that, with probability tending to 1,*

$$\sup_{\lambda \in \Lambda_n} |L(\lambda) - \widehat{L}(\lambda) - \delta_n| = O_P\left(\frac{k_n}{n^{1-c_2}}\right) + O_P\left(\frac{k_n}{\sqrt{n}}\right). \tag{9}$$

**4. Multi-stage methods.** The multi-stage methods use the following steps. As mentioned earlier, we randomly split the data into three parts, $\mathcal{D}_1$, $\mathcal{D}_2$ and $\mathcal{D}_3$, which we take to be of equal size:

1. Stage I. Use $\mathcal{D}_1$ to find $\widehat{S}_n(\lambda)$, for each $\lambda$.
2. Stage II. Use $\mathcal{D}_2$ to find $\widehat{\lambda}$ by cross-validation, and let $\widehat{S}_n = \widehat{S}_n(\widehat{\lambda})$.
3. Stage III. Use $\mathcal{D}_3$ to find the least-squares estimate $\widehat{\beta}$ for the model $\widehat{S}_n$.
   Let

$$\widehat{D}_n = \{j \in \widehat{S}_n : |T_j| > c_n\},$$

where $T_j$ is the usual $t$-statistic, $c_n = z_{\alpha/2m}$ and $m = |\widehat{S}_n|$.

4.1. *The lasso.* The lasso estimator [Tibshirani (1996)] $\widetilde{\beta}(\lambda)$ minimizes

$$M_\lambda(\lambda) = \sum_{i=1}^n (Y_i - X_i^T \beta)^2 + \lambda \sum_{j=1}^p |\beta_j|$$

and let $\widehat{S}_n(\lambda) = \{j : \widetilde{\beta}_j(\lambda) \neq 0\}$. Recall that $\widehat{\beta}(\lambda)$ is the least-squares estimator using the covariates in $\widehat{S}_n(\lambda)$.

Let $k_n = A \log n$ where $A > 0$ is a positive constant.

THEOREM 4.1. *Assume that* (A1)–(A5) *hold. Let* $\Lambda_n = \{\lambda : |\widehat{S}_n(\lambda)| \leq k_n\}$. *Then:*

1. *The true loss overfits:* $\mathbb{P}(D \subset \widehat{S}_n(\lambda_*)) \to 1$ *where* $\lambda_* = \arg\min_{\lambda \in \Lambda_n} L(\lambda)$.

2. *Cross-validation also overfits:* $\mathbb{P}(D \subset \widehat{S}_n(\widehat{\lambda})) \to 1$ *where* $\widehat{\lambda} = \arg\min_{\lambda \in \Lambda_n} \widehat{L}(\lambda)$.
3. *Type* I *error is controlled:* $\limsup_{n\to\infty} \mathbb{P}(D^c \cap \widehat{D}_n \neq \varnothing) \leq \alpha$.

If we let $\alpha = \alpha_n \to 0$, then $\widehat{D}_n$ is consistent for variable selection.

THEOREM 4.2. *Assume that* (A1)–(A5) *hold. Let* $\alpha_n \to 0$ *and* $\sqrt{n}\alpha_n \to \infty$. *Then, the multi-stage lasso is consistent,*

$$\mathbb{P}(\widehat{D}_n = D) \to 1. \tag{10}$$

The next result follows directly. The proof is thus omitted.

THEOREM 4.3. *Assume that* (A1)–(A5) *hold. Let* $\alpha$ *be fixed. Then,* $(\widehat{D}_n, \widehat{S}_n)$ *forms a confidence sandwich*

$$\liminf_{n\to\infty} \mathbb{P}(\widehat{D}_n \subset D \subset \widehat{S}_n) \geq 1 - \alpha. \tag{11}$$

REMARK 4.4. This confidence sandwich is expected to be conservative in the sense that the coverage can be much larger than $1 - \alpha$.

4.2. *Stepwise regression.* Let $k_n = A \log n$ for some $A > 0$. The version of stepwise regression we consider is as follows:

1. Initialize: $\text{Res} = Y$, $\lambda = 0$, $\widehat{Y} = 0$ and $\widehat{S}_n(\lambda) = \varnothing$.
2. Let $\lambda \leftarrow \lambda + 1$. Compute $\widehat{\mu}_j = n^{-1}\langle X_j, \text{Res}\rangle$ for $j = 1, \ldots, p$.
3. Let $J = \arg\max_j |\widehat{\mu}_j|$. Set $\widehat{S}_n(\lambda) = \{\widehat{S}_n(\lambda - 1), J\}$. Set $\widehat{Y} = X_\lambda \widehat{\beta}(\lambda)$ where $\widehat{\beta}_\lambda = (X_\lambda^T X_\lambda)^{-1} X_\lambda^T Y$, and let $\text{Res} = Y - \widehat{Y}$.
4. If $\lambda = k_n$, stop. Otherwise, go to step 2.

For technical reasons, we assume that the final estimator $x^T \widehat{\beta}$ is truncated to be no larger than $B$. Note that $\lambda$ is discrete and $\Lambda_n = \{0, 1, \ldots, k_n\}$.

THEOREM 4.5. *With* $\widehat{S}_n(\lambda)$ *defined as above, the statements of Theorems 4.1, 4.2 and 4.3 hold.*

4.3. *Marginal regression.* This is probably the oldest, simplest and most common method. It is quite popular in gene expression analysis. It used to be regarded with some derision but has enjoyed a revival. A version appears in a recent paper by Fan and Lv (2008). Let $\widehat{S}_n(\lambda) = \{j : |\widehat{\mu}_j| \geq \lambda\}$ where $\widehat{\mu}_j = n^{-1}\langle Y, X_{\bullet j}\rangle$.

Let $\mu_j = \mathbb{E}(\widehat{\mu}_j)$, and let $\mu_{(j)}$ denote the value of $\mu$ ordered by their absolute values,

$$|\mu_{(1)}| \geq |\mu_{(2)}| \geq \cdots.$$

THEOREM 4.6. *Let $k_n \to \infty$ with $k_n = o(\sqrt{n})$. Let $\Lambda_n = \{\lambda : |\widehat{S}_n(\lambda)| \leq k_n\}$. Assume that*

$$\min_{j \in D} |\mu_j| > |\mu_{(k_n)}|. \tag{12}$$

*Then, the statements of Theorems 4.1, 4.2 and 4.3 hold.*

Assumption (12) limits the degree of unfaithfulness (small partial correlations induced by cancellation of parameters). Large values of $k_n$ weaken assumption (12), thus making the method more robust to unfaithfulness, but at the expense of lower power. Fan and Lv (2008) make similar assumptions. They assume that there is a $C > 0$ such that $|\mu_j| \geq C|\beta_j|$ for all $j$, which rules out unfaithfulness. However, they do not explicitly relate the values of $\mu_j$ for $j \in D$ to the values outside $D$ as we have done. On the other hand, they assume that $Z = \Sigma^{-1/2} X$ has a spherically symmetric distribution. Under this assumption and their faithfulness assumption, they deduce that the $\mu_j$'s outside $D$ cannot strongly dominate the $\mu_j$'s within $D$. We prefer to simply make this an explicit assumption without placing distributional assumptions on $X$. At any rate, any method that uses marginal regressions as a starting point must make some sort of faithfulness assumptions to succeed.

4.4. *Modifications.* Let us now discuss a few modifications of the basic method. First, consider splitting the data only into two groups, $\mathcal{D}_1$ and $\mathcal{D}_2$. Then do the following steps:

1. Stage I. Find $\widehat{S}_n(\lambda)$ for $\lambda \in \Lambda_n$, where $|\widehat{S}_n(\lambda)| \leq k_n$ for each $\lambda \in \Lambda_n$ using $\mathcal{D}_1$.
2. Stage II. Find $\widehat{\lambda}$ by cross-validation, and let $\widehat{S}_n = \widehat{S}_n(\widehat{\lambda})$ using $\mathcal{D}_2$.
3. Stage III. Find the least-squares estimate $\widehat{\beta}_{\widehat{S}_n}$ using $\mathcal{D}_2$. Let $\widehat{D}_n = \{j \in \widehat{S}_n : |T_j| > c_n\}$, where $T_j$ is the usual $t$-statistic.

THEOREM 4.7. *Choosing*

$$c_n = \frac{\log \log n \sqrt{2k_n \log(2p_n)}}{\alpha} \tag{13}$$

*controls asymptotic type* I *error.*

The critical value in (13) is hopelessly large and it does not appear it can be substantially reduced. We present this mainly to show the value of the extra data-splitting step. It is tempting to use the same critical value as in the tri-split case, namely $c_n = z_{\alpha/2m}$ where $m = |\widehat{S}_n|$, but we suspect this will not work in general. However, it may work under extra conditions.

**5. Application.** As an example, we illustrate an analysis based on part of the osteoporotic fractures in men study [MrOS, Orwoll et al. (2005)]. A sample of 860 men were measured at a large number of genes and outcome measures. We consider only 296 SNPs which span 30 candidate genes for bone mineral density. An aim of the study was to identify genes associated with bone mineral density that could help in understanding the genetic basis of osteoporosis in men. Initial analyses of this subset of the data revealed no SNPs with a clear pattern of association with the phenotype; however, three SNPs, numbered (67), (277) and (289), exhibited some association in the screening of the data. To further explore the effacacy of the lasso screen and clean procedure, we modified the phenotype to enhance this weak signal and then reanalyzed the data to see if we could detect this planted signal.

We were interested in testing for main effects and pairwise interactions in these data; however, including all interactions results in a model with 43,660 additional terms, which is not practical for this sample size. As a compromise, we selected 2 SNPs per gene to model potential interaction effects. This resulted in a model with a total of 2066 potential coefficients, including 296 main effects and 1770 interaction terms. With this model, our initial screen detected 10 terms, including the 3 enhanced signals, 2 other main effects and 5 interactions. After cleaning, the final model detected the 3 enhanced signals and no other terms.

**6. Simulations.** To further explore the screen and clean procedures, we conducted simulation experiments with four models. For each model $Y_i = X_i^T \beta + \varepsilon_i$ where the measurement errors, $\varepsilon_i$ and $\varepsilon_{ij}^*$, are i.i.d. normal$(0,1)$ and the covariates $X_{ij}$'s are normal$(0,1)$ (except for model D). Models differ in how $Y_i$ is linked to $X_i$ and the dependence structure of the $X_i$'s. Models A, B and C explore scenarios with moderate and large $p$, while model D focuses on confounding and unfaithfullness, as follows:

(A) Null model: $\beta = (0, \ldots, 0)$ and the $X_{ij}$'s are i.i.d.
(B) Triangle model: $\beta_j = \delta(10-j), j = 1, \ldots, 10, \beta_j = 0, j > 10$ and $X_{ij}$'s are i.i.d.
(C) Correlated Triangle model: as B, but with $X_{ij(+1)} = \rho X_{ij} + (1-\rho^2)^{1/2} \varepsilon_{ij}^*$, for $j > 1$, and $\rho = 0.5$.
(D) Unfaithful model: $Y_i = \beta_1 X_{i1} + \beta_2 X_{i2} + \varepsilon_i$, for $\beta_1 = -\beta_2 = 10$, where the $X_{ij}$'s are i.i.d. for $j = \{1, 5, 6, 7, 8, 9, 10\}$, but $X_{i2} = \rho X_{i1} + \tau \varepsilon_{i2}^*$, $X_{i3} = \rho X_{i1} + \tau \varepsilon_{i10}^*$, and $X_{i4} = \rho X_{i2} + \tau \varepsilon_{i11}^*$, for $\tau = 0.01$ and $\rho = 0.95$.

We used a maximum model size of $k_n = n^{1/2}$ which technically goes beyond the theory but works well in practice. Prior to analysis, the covariates are scaled so that each has mean 0 and variance 1. The tests were initially performed using a third of the data for each of the 3 stages of the procedure (Table 1, top half, 3 splits). For models A, B and C, each approach has type I error less than $\alpha$, except the stepwise procedure which has trouble with model C when $n = p = 100$. We also calculated the false positive rate and found it to be very low (about $10^{-4}$ when $p = 100$ and $10^{-5}$ when $p = 1000$) indicating that even when a type I error occurs, only a very small number of terms are included erroneously. The lasso screening procedure exhibited a slight power advantage over the stepwise procedure. Both methods dominated the marginal approach. The Markov dependence structure in model C clearly challenged the marginal approach. For model D, none of the approaches controlled the type I error.

TABLE 1
*Size and power of screen and clean procedures using lasso, stepwise and marginal regression for the screening step. For all procedures $\alpha = 0.05$. For $p = 100$, $\delta = 0.5$ and for $p = 1000$, $\delta = 1.5$. Reported power is $\pi_{\mathrm{av}}$. The top 8 rows of simulations were conducted using three stages as described in Section 4, with a third of the data used for each stage. The bottom 8 rows of simulations were conducted splitting the data in half, using the first portion with leave-one-out cross validation for stages 1 and 2 and the second portion for cleaning*

| | | | | Size | | | Power | | |
|---|---|---|---|---|---|---|---|---|---|
| Splits | $n$ | $p$ | Model | Lasso | Step | Marg | Lasso | Step | Marg |
| 2 | 100 | 100 | A | 0.005 | 0.001 | 0.004 | 0 | 0 | 0 |
| 2 | 100 | 100 | B | 0.01 | 0.02 | 0.03 | 0.62 | 0.62 | 0.31 |
| 2 | 100 | 100 | C | 0.001 | 0.01 | 0.01 | 0.77 | 0.57 | 0.21 |
| 2 | 100 | 10 | D | 0.291 | 0.283 | 0.143 | 0.08 | 0.08 | 0.04 |
| 2 | 100 | 1000 | A | 0.001 | 0.002 | 0.010 | 0 | 0 | 0 |
| 2 | 100 | 1000 | B | 0.002 | 0.020 | 0.010 | 0.17 | 0.09 | 0.11 |
| 2 | 100 | 1000 | C | 0.02 | 0.14 | 0.01 | 0.27 | 0.15 | 0.11 |
| 2 | 1000 | 10 | D | 0.291 | 0.283 | 0.143 | 0.08 | 0.08 | 0.04 |
| 3 | 100 | 100 | A | 0.040 | 0.050 | 0.030 | 0 | 0 | 0 |
| 3 | 100 | 100 | B | 0.02 | 0.01 | 0.02 | 0.91 | 0.90 | 0.56 |
| 3 | 100 | 100 | C | 0.03 | 0.04 | 0.03 | 0.91 | 0.88 | 0.41 |
| 3 | 100 | 10 | D | 0.382 | 0.343 | 0.183 | 0.16 | 0.18 | 0.09 |
| 3 | 100 | 1000 | A | 0.035 | 0.045 | 0.040 | 0 | 0 | 0 |
| 3 | 100 | 1000 | B | 0.045 | 0.020 | 0.035 | 0.57 | 0.66 | 0.29 |
| 3 | 100 | 1000 | C | 0.06 | 0.070 | 0.020 | 0.74 | 0.65 | 0.19 |
| 3 | 1000 | 10 | D | 0.481 | 0.486 | 0.187 | 0.17 | 0.17 | 0.13 |

To determine the sensitivity of the approach to using distinct data for each stage of the analysis, simulations were conducted screening on the first half of the data and cleaning on the second half (2 splits). The tuning parameter was selected using leave-one-out cross validation (Table 1, bottom half). As expected, this approach lead to a dramatic increase in the power of all the procedures. More surprising is the fact that the type I error was near $\alpha$ or below for models A, B and C. Clearly this approach has advantages over data splitting and merits further investigation.

A natural competitor to screen and clean procedure is a two-stage adaptive lasso [Zou (2006)]. In our implementation, we split the data and used one half for each stage of the analysis. At stage one, leave-one-out cross validation lasso screens the data. In stage two, the adaptive lasso, with weights $w_j = |\widehat{\beta}_j|^{-1}$, cleans the data. The tuning parameter for the lasso was again chosen using leave-one-out cross validation. Table 2 provides the size, power and false positive rate (FPR) for this procedure. Naturally, the adaptive lasso does not control the size of the test, but the FPR is small. The power of the test is greater than we found for our lasso screen and clean procedure, but this extra power comes at the cost of a much higher type I error rate.

**7. Proofs.** Recall that if $A$ is a square matrix, then $\phi(A)$ and $\Phi(A)$ denote the smallest and largest eigenvalues of $A$. Throughout the proofs we make use of the following fact. If $v$ is a vector and $A$ is a square matrix, then

$$(14) \qquad \phi(A)\|v\|^2 \leq v^T A v \leq \Phi(A)\|v\|^2.$$

We use the following standard tail bound: if $Z \sim N(0,1)$, then $\mathbb{P}(|Z| > t) \leq t^{-1} e^{-t^2/2}$. We will also use the following results about the lasso from Meinshausen and Yu (2008). Their results are stated and proved for fixed $X$ but, under the conditions (A1)–(A5), it is easy to see that their conditions hold with probability tending to one and so their results hold for random $X$ as well.

THEOREM 7.1 [Meinshausen and Yu (2008)]. *Let $\widetilde{\beta}(\lambda)$ be the lasso estimator.*

Table 2
*Size, power and false positive rate (FPR) of two-stage adaptive lasso procedure*

| $n$ | $p$ | Model | Size | Power | FPR |
|---|---|---|---|---|---|
| 100 | 100 | A | 0.93 | 0 | 0.032 |
| 100 | 100 | B | 0.84 | 0.97 | 0.034 |
| 100 | 100 | C | 0.81 | 0.96 | 0.031 |
| 100 | 10 | D | 0.67 | 0.21 | 0.114 |
| 100 | 1000 | A | 0.96 | 0 | 0.004 |
| 100 | 1000 | B | 0.89 | 0.65 | 0.004 |
| 100 | 1000 | C | 0.76 | 0.77 | 0.002 |
| 1000 | 10 | D | 0.73 | 0.24 | 0.013 |

1. *The squared error satisfies*

(15) $$\mathbb{P}\bigg(\|\widetilde{\beta}(\lambda) - \beta\|_2^2 \leq \frac{2\lambda^2 s}{n^2\kappa^2} + \frac{cm\log p_n}{n\phi_n^2(m)}\bigg) \to 1,$$

*where $m = |\widehat{S}_n(\lambda)|$ and $c > 0$ is a constant.*

2. *The size of $\widehat{S}_n(\lambda)$ satisfies*

(16) $$\mathbb{P}\bigg(|\widehat{S}_n(\lambda)| \leq \frac{\tau^2 C n^2}{\lambda^2}\bigg) \to 1,$$

*where $\tau^2 = \mathbb{E}(Y_i^2)$.*

PROOF OF LEMMA 3.1. Let $D \subset M$ and $\phi = \phi(n^{-1}X_M^T X_M)$. Then,

$$L(\widehat{\beta}_M) = \frac{1}{n}\varepsilon^T X_M (X_M^T X_M)^{-1} X_M^T \varepsilon \leq \frac{1}{n^2\phi}\|X_M^T\varepsilon\|^2 = \frac{1}{n\phi}\sum_{j\in M} Z_j^2,$$

where $Z_j = n^{-1/2} X_{\bullet j}^T \varepsilon$. Conditional on $X$, $Z_i \sim N(0, a_j^2)$ where $a_j^2 = n^{-1} \times \sum_{i=1}^n X_{ij}^2$. Let $A_n^2 = \max_{1\leq j \leq p_n} a_j^2$. By Hoeffding's inequality, (A2) and (A5), $\mathbb{P}(E_n) \to 1$ where $E_n = \{A_n \leq \sqrt{2}\}$. So

$$\mathbb{P}\bigg(\max_{1\leq j\leq p_n} |Z_j| > \sqrt{4\log p_n}\bigg)$$

$$= \mathbb{P}\bigg(\max_{1\leq j\leq p_n} |Z_j| > \sqrt{4\log p_n}, E_n\bigg) + \mathbb{P}\bigg(\max_{1\leq j\leq p_n} |Z_j| > \sqrt{4\log p_n}, E_n^c\bigg)$$

$$\leq \mathbb{P}\bigg(\max_{1\leq j\leq p_n} |Z_j| > \sqrt{4\log p_n}, E_n\bigg) + \mathbb{P}(E_n^c)$$

$$\leq \mathbb{P}\bigg(A_n \max_{1\leq j\leq p_n} \frac{|Z_j|}{a_j} > \sqrt{4\log p_n}, E_n\bigg) + o(1)$$

$$\leq \mathbb{P}\bigg(\max_{1\leq j\leq p_n} \frac{|Z_j|}{a_j} > \sqrt{2\log p_n}\bigg) + o(1)$$

$$= \mathbb{E}\bigg(\mathbb{P}\bigg(\max_{1\leq j\leq p_n} \frac{|Z_j|}{a_j} > \sqrt{2\log p_n}\bigg)\bigg|X\bigg) + o(1)$$

$$\leq O\bigg(\frac{1}{\sqrt{2\log p_n}}\bigg) + o(1) = o(1).$$

But $\sum_{j\in M} Z_j^2 \leq m \max_{1\leq j\leq p_n} Z_j^2$ and (6) follows.

Now we lower bound $L(\widehat{\beta}_M)$. Let $M$ be such that $D \not\subset M$. Let $A = \{j : \widehat{\beta}_M(j) \neq 0\} \cup D$. Then, $|A| \leq m + s$. Therefore, with probability tending to 1,

$$L(\widehat{\beta}_M) = \frac{1}{n}(\widehat{\beta}_M - \beta)^T X^T X (\widehat{\beta}_M - \beta) = \frac{1}{n}(\widehat{\beta}_M - \beta)^T X_A^T X_A (\widehat{\beta}_M - \beta)$$

$$\geq \phi_n(m+s)\|\widehat{\beta}_M - \beta\|^2 = \phi_n(m+s)\sum_{j\in A}(\widehat{\beta}_M(j) - \beta(j))^2$$

$$\geq \phi_n(m+s)\sum_{j\in D\cap M^c}(0 - \beta(j))^2 \geq \phi_n(m+s)\psi^2.$$
□

PROOF OF THEOREM 3.2. Let $\widetilde{Y}$ denote the responses, and $\widetilde{X}$ the design matrix, for the second half of the data. Then, $\widetilde{Y} = \widetilde{X}\beta + \widetilde{\varepsilon}$. Now

$$L(\lambda) = \frac{1}{n}(\widehat{\beta}(\lambda) - \beta)^T X^T X(\widehat{\beta}(\lambda) - \beta) = (\widehat{\beta}(\lambda) - \beta)^T \widehat{\Sigma}_n(\widehat{\beta}(\lambda) - \beta)$$

and

$$\widehat{L}(\lambda) = n^{-1}\|\widetilde{Y} - \widetilde{X}\widehat{\beta}(\lambda)\|^2 = (\widehat{\beta}(\lambda) - \beta)^T \widetilde{\Sigma}_n(\widehat{\beta}(\lambda) - \beta) + \delta_n + \frac{2}{n}\langle \widetilde{\varepsilon}, \widetilde{X}(\widehat{\beta}(\lambda) - \beta)\rangle,$$

where $\delta_n = \|\widetilde{\varepsilon}\|^2/n$, and $\widehat{\Sigma}_n = n_1^{-1}X^T X$ and $\widetilde{\Sigma}_n = n^{-1}\widetilde{X}^T\widetilde{X}$. By Hoeffding's inequality

$$\mathbb{P}(|\widehat{\Sigma}_n(j,k) - \widetilde{\Sigma}_n(j,k)| > \varepsilon) \leq e^{-nc\varepsilon^2}$$

for some $c > 0$, and so

$$\mathbb{P}\left(\max_{jk}|\widehat{\Sigma}_n(j,k) - \widetilde{\Sigma}_n(j,k)| > \varepsilon\right) \leq p_n^2 e^{-nc\varepsilon^2}.$$

Choose $\varepsilon_n = 4/(cn^{1-c_2})$. It follows that

$$\mathbb{P}\left(\max_{jk}|\widehat{\Sigma}_n(j,k) - \widetilde{\Sigma}_n(j,k)| > \frac{4}{cn^{1-c_2}}\right) \leq e^{-2n^{c_2}} \to 0.$$

Note that

$$|\{j : \widehat{\beta}_j(\lambda) \neq 0\} \cup \{j : \beta_j \neq 0\}| \leq k_n + s.$$

Hence, with probability tending to 1

$$|L(\lambda) - \widehat{L}(\lambda) - \delta_n| \leq \frac{4}{cn^{1-c_2}}\|\widehat{\beta}(\lambda) - \beta\|_1^2 + 2\xi_n(\lambda)$$

for all $\lambda \in \Lambda_n$, where

$$\xi_n(\lambda) = \frac{1}{n}\sum_{i\in I_2}\widetilde{\varepsilon}_i\mu_i(\lambda)$$

and $\mu_i(\lambda) = \widetilde{X}_i^T(\widehat{\beta}(\lambda) - \beta)$. Now $\|\widehat{\beta}(\lambda) - \beta\|_1^2 = O_P((k_n + s)^2)$ since $\|\widehat{\beta}(\lambda)\|^2 = O_P(k_n/\phi(k_n))$. Thus, $\|\widehat{\beta}(\lambda) - \beta\|_1 \leq C(k_n + s)$ with probability tending to 1, for some $C > 0$. Also, $|\mu_i(\lambda)| \leq B\|\widehat{\beta}(\lambda) - \beta\|_1 \leq BC(k_n + s)$ with probability tending to 1. Let $W \sim N(0,1)$. Conditional on $\mathcal{D}_1$,

$$|\xi_n(\lambda)| \stackrel{d}{=} \frac{\sigma}{\sqrt{n}}\sqrt{\sum_{i=1}^n \mu_i^2(\lambda)}|W| \leq \frac{\sigma}{\sqrt{n}}BC(k_n + s)|W|,$$

so $\sup_{\lambda \in \Lambda_n} |\xi_n(\lambda)| = O_P(k_n/\sqrt{n})$. □

PROOF OF THEOREM 4.1. 1. Let $\lambda_n = \tau n \sqrt{C/k_n}$, $M = \widehat{S}_n(\lambda_n)$ and $m = |M|$. Then, $\mathbb{P}(m \leq k_n) \to 1$ due to (16). Hence, $\mathbb{P}(\lambda_n \in \Lambda_n) \to 1$. From (15),

$$\|\widetilde{\beta}(\lambda_n) - \beta\|_2^2 \leq O\left(\frac{1}{k_n}\right) + O_P\left(\frac{k_n \log p_n}{n}\right) = o_P(1).$$

Hence, $\|\widetilde{\beta}(\lambda_n) - \beta\|_\infty^2 = o_P(1)$. So, for each $j \in D$,

$$|\widetilde{\beta}_j(\lambda_n)| \geq |\beta_j| - |\widetilde{\beta}_j(\lambda_n) - \beta_j| \geq \psi + o_P(1)$$

and hence, $\mathbb{P}(\min_{j \in D} |\widetilde{\beta}_j(\lambda_n)| > 0) \to 1$. Therefore, $\Gamma_n = \{\lambda \in \Lambda_n : D \subset \widehat{S}_n(\lambda)\}$ is nonempty. By Lemma 3.1,

(17) $$L(\lambda_n) \leq cm \log p_n / (n\phi(m)) = O_P(k_n \log p_n / n).$$

On the other hand, from Lemma 3.1,

(18) $$\mathbb{P}\left(\inf_{\lambda \in \Lambda_n \cap \Gamma_n^c} L(\widehat{\beta}_\lambda) > \psi^2 \phi(k_n)\right) \to 1.$$

Now, $n\phi_n(k_n)/(k_n \log p_n) \to \infty$, and so (17) and (18) imply that

$$\mathbb{P}\left(\inf_{\lambda \in \Lambda_n \cap \Gamma_n^c} L(\widehat{\beta}_\lambda) > L(\lambda_n)\right) \to 1.$$

Thus, if $\lambda_*$ denotes the minimizer of $L(\lambda)$ over $\Lambda_n$, we conclude that $\mathbb{P}(\lambda_* \in \Gamma_n) \to 1$, and hence, $\mathbb{P}(D \subset \widehat{S}_n(\lambda_*)) \to 1$.

2. This follows from part 1 and Theorem 3.2.

3. Let $A = \widehat{S}_n \cap D^c$. We want to show that

$$\mathbb{P}\left(\max_{j \in A} |T_j| > c_n\right) \leq \alpha + o(1).$$

Now,

$$\mathbb{P}\left(\max_{j \in A} |T_j| > c_n\right) = \mathbb{P}\left(\max_{j \in A} |T_j| > c_n, D \subset \widehat{S}_n\right) + \mathbb{P}\left(\max_{j \in A} |T_j| > c_n, D \not\subset \widehat{S}_n\right)$$

$$\leq \mathbb{P}\left(\max_{j \in A} |T_j| > c_n, D \subset \widehat{S}_n\right) + \mathbb{P}(D \not\subset \widehat{S}_n)$$

$$\leq \mathbb{P}\left(\max_{j \in A} |T_j| > c_n, D \subset \widehat{S}_n\right) + o(1).$$

Conditional on $(\mathcal{D}_1, \mathcal{D}_2)$, $\widehat{\beta}_A$ is normally distributed with mean 0 and variance matrix $\sigma^2 (X_A^T X_A)^{-1}$ when $D \subset \widehat{S}_n$. Recall that

$$T_j(M) = \frac{e_j^T (X_M^T X_M)^{-1} X_M^T Y}{\widehat{\sigma} \sqrt{e_j^T (X_M^T X_M)^{-1} e_j}} = \frac{\widehat{\beta}_{M,j}}{s_j},$$

where $M = \widehat{S}_n$, $s_j^2 = \widehat{\sigma}^2 e_j^T (X_M^T X_M)^{-1} e_j$ and $e_j = (0, \ldots, 0, 1, 0, \ldots, 0)^T$, where the 1 is in the $j$th coordinate. When $D \subset \widehat{S}_n$, each $T_j$, for $j \in A$, has a $t$-distribution with $n - m$ degrees of freedom where $m = |\widehat{S}_n|$. Also, $c_n/t_{\alpha/2m} \to 1$ where $t_u$ denotes the upper tail critical value for the $t$-distribution. Hence,

$$\mathbb{P}\left(\max_{j \in A} |T_j| > c_n, D \subset \widehat{S}_n | \mathcal{D}_1, \mathcal{D}_2\right)$$
$$= \mathbb{P}\left(\max_{j \in A} |T_j| > t_{\alpha/2m}, D \subset \widehat{S}_n | \mathcal{D}_1, \mathcal{D}_2\right) + a_n$$
$$\leq \alpha + a_n,$$

where $a_n = o(1)$, since $|A| \leq m$. It follows that

$$\mathbb{P}\left(\max_{j \in A} |T_j| > c_n, D \subset \widehat{S}_n\right) \leq \alpha + o(1). \qquad \square$$

PROOF OF THEOREM 4.2.  From Theorem 4.1, $\mathbb{P}(\widehat{D}_n \cap D^c \neq \varnothing) \leq \alpha_n$ and so $\mathbb{P}(\widehat{D}_n \cap D^c \neq \varnothing) \to 0$. Hence, $\mathbb{P}(\widehat{D}_n \subset D) \to 1$. It remains to be shown that

(19) $$\mathbb{P}(D \subset \widehat{D}_n) \to 1.$$

The test statistic for testing $\beta_j = 0$ when $\widehat{S}_n = M$ is

$$T_j(M) = \frac{e_j^T (X_M^T X_M)^{-1} X_M^T Y}{\widehat{\sigma} \sqrt{e_j^T (X_M^T X_M)^{-1} e_j}}.$$

For simplicity in the proof, let us take $\widehat{\sigma} = \sigma$, the extension to unknown $\sigma$ being straightforward. Let $j \in D$, $\mathcal{M} = \{M : |M| \leq k_n, D \subset M\}$. Then,

$$\mathbb{P}(j \notin \widehat{D}_n) = \mathbb{P}(j \notin \widehat{D}_n, D \subset \widehat{S}_n) + \mathbb{P}(j \notin \widehat{D}_n, D \not\subset \widehat{S}_n)$$
$$\leq \mathbb{P}(j \notin \widehat{D}_n, D \subset \widehat{S}_n) + \mathbb{P}(D \not\subset \widehat{S}_n)$$
$$= \mathbb{P}(j \notin \widehat{D}_n, D \subset \widehat{S}_n) + o(1)$$
$$= \sum_{M \in \mathcal{M}} \mathbb{P}(j \notin \widehat{D}_n, \widehat{S}_n = M) + o(1)$$
$$\leq \sum_{M \in \mathcal{M}} \mathbb{P}(|T_j(M)| < c_n, \widehat{S}_n = M) + o(1)$$
$$\leq \sum_{M \in \mathcal{M}} \mathbb{P}(|T_j(M)| < c_n) + o(1).$$

Conditional on $\mathcal{D}_1 \cup \mathcal{D}_2$, for each $M \in \mathcal{M}$, $T_j(M) = (\beta_j/s_j) + Z$, where $Z \sim N(0,1)$. Without loss of generality, assume that $\beta_j > 0$. Hence,

$$\mathbb{P}(|T_j(M)| < c_n | \mathcal{D}_1 \cup \mathcal{D}_2) = \mathbb{P}\left(-c_n - \frac{\beta_j}{s_j} < Z < c_n - \frac{\beta_j}{s_j}\right).$$

Fix a small $\varepsilon > 0$. Note that $s_j^2 \leq \sigma^2/(n\kappa)$. It follows that, for all large $n$, $c_n - \beta_j/s_j < -\varepsilon\sqrt{n}$. So,

$$\mathbb{P}(|T_j(M)| < c_n | \mathcal{D}_1 \cup \mathcal{D}_2) \leq \mathbb{P}(Z < -\varepsilon\sqrt{n}) \leq e^{-n\varepsilon^2/2}.$$

The number of models in $\mathcal{M}$ is

$$\sum_{j=0}^{k_n} \binom{p_n - s}{j - s} \leq k_n \binom{p_n - s}{k_n - s} \leq k_n \left(\frac{(p_n - s)e}{k_n - s}\right)^{k_n - s} \leq k_n p_n^{k_n},$$

where we used the inequality

$$\binom{n}{k} \leq \left(\frac{ne}{k}\right)^k.$$

So,

$$\sum_{M \in \mathcal{M}} \mathbb{P}(|T_j(M)| < c_n | \mathcal{D}_1 \cup \mathcal{D}_2) \leq k_n p_n^{k_n} e^{-n\varepsilon^2} \to 0$$

by (A2). We have thus shown that $\mathbb{P}(j \notin \widehat{D}_n) \to 0$, for each $j \in D$. Since $|D|$ is finite, it follows that $\mathbb{P}(j \notin \widehat{D}_n \text{ for some } j \in D) \to 0$ and hence (19). $\square$

PROOF OF THEOREM 4.5. A simple modification of Theorem 3.1 of Barron et al. (2008) shows that

$$L(k_n) = \frac{1}{n}\|\widehat{Y}_{k_n} - X\beta\|^2 = o_P(1).$$

[The modification is needed because Barron et al. (2008) require $Y$ to be bounded while we have assumed that $Y$ is normal. By a truncation argument, we can still derive the bound on $L(k_n)$.] So

$$\|\widehat{\beta}_{k_n} - \beta\|^2 \leq \frac{L(k_n)}{\phi_n(k_n + s)} \leq \frac{L(k_n)}{\kappa} = o_P(1).$$

Hence, for any $\varepsilon > 0$, with probability tending to 1, $\|\widehat{\beta}(k_n) - \beta\|^2 < \varepsilon$ so that $|\widehat{\beta}_j| > \psi/2 > 0$, for all $j \in D$. Thus, $\mathbb{P}(D \subset \widehat{S}_n(k_n)) \to 1$. The remainder of the proof of part 1 is the same as in Theorem 4.1. Part 2 follows from the previous result together with Theorem 3.2. The proof of part 3 is the same as for Theorem 4.1. $\square$

PROOF OF THEOREM 4.6. Note that $\widehat{\mu}_j - \mu_j = n^{-1}\sum_{i=1}^n X_{ij}\varepsilon_i$. Hence, $\widehat{\mu}_j - \mu_j \sim N(0, 1/n)$. So, for any $\delta > 0$,

$$\mathbb{P}\left(\max_j |\widehat{\mu}_j - \mu_j| > \delta\right) \leq \sum_{j=1}^{p_n} \mathbb{P}(|\widehat{\mu}_j - \mu_j| > \delta)$$

$$\leq \frac{p_n}{\delta\sqrt{n}} e^{-n\delta^2/2} \leq \frac{c_1 e^{n^{c_2}}}{\delta\sqrt{n}} e^{-n\delta^2/2} \to 0.$$

By (12), conclude that $D \subset \widehat{S}_n(\lambda)$ when $\lambda = \widehat{\mu}_{(k_n)}$. The remainder of the proof is the same as the proof of Theorem 4.5. $\square$

PROOF OF THEOREM 4.7. Let $A = \widehat{S}_n \cap D^c$. We want to show that

$$\mathbb{P}\left(\max_{j \in A} |T_j| > c_n\right) \leq \alpha + o(1).$$

For fixed $A$, $\widehat{\beta}_A$ is normal with mean 0 but this is not true for random $A$. Instead we need to bound $T_j$. Recall that

$$T_j(M) = \frac{e_j^T(X_M^T X_M)^{-1} X_M^T Y}{\widehat{\sigma}\sqrt{e_j^T(X_M^T X_M)^{-1} e_j}} = \frac{\widehat{\beta}_{M,j}}{s_j},$$

where $M = \widehat{S}_n$, $s_j^2 = \widehat{\sigma}^2 e_j^T(X_M^T X_M)^{-1} e_j$ and $e_j = (0, \ldots, 0, 1, 0, \ldots, 0)^T$ where the 1 is in the $j$th coordinate. The probabilities that follow are conditional on $\mathcal{D}_1$ but this is supressed for notational convenience. First, write

$$\mathbb{P}\left(\max_{j \in A} |T_j| > c_n\right) = \mathbb{P}\left(\max_{j \in A} |T_j| > c_n, D \subset \widehat{S}_n\right) + \mathbb{P}\left(\max_{j \in A} |T_j| > c_n, D \not\subset \widehat{S}_n\right)$$

$$\leq \mathbb{P}\left(\max_{j \in A} |T_j| > c_n, D \subset \widehat{S}_n\right) + \mathbb{P}(D \not\subset \widehat{S}_n)$$

$$\leq \mathbb{P}\left(\max_{j \in A} |T_j| > c_n, D \subset \widehat{S}_n\right) + o(1).$$

When $D \subset \widehat{S}_n$,

$$\widehat{\beta}_{\widehat{S}_n} = \beta_{\widehat{S}_n} + \left(\frac{1}{n} X_{\widehat{S}_n}^T X_{\widehat{S}_n}\right)^{-1} \frac{1}{n} X_{\widehat{S}_n}^T \varepsilon = \beta_{\widehat{S}_n} + Q_{\widehat{S}_n} \gamma_{\widehat{S}_n},$$

where $Q_{\widehat{S}_n} = ((1/n) X_{\widehat{S}_n}^T X_{\widehat{S}_n})^{-1}$, $\gamma_{\widehat{S}_n} = n^{-1} X_{\widehat{S}_n}^T \varepsilon$, and $\beta_{\widehat{S}_n}(j) = 0$, for $j \in A$. Now, $s_j^2 \geq \widehat{\sigma}^2/(nC)$ so that

$$|T_j| \leq \frac{\sqrt{nC}|\widehat{\beta}_{\widehat{S}_n,j}|}{\widehat{\sigma}} \leq \frac{\sqrt{n \log \log n}|\widehat{\beta}_{\widehat{S}_n,j}|}{\widehat{\sigma}}$$

for $j \in \widehat{S}_n$. Therefore,

$$\mathbb{P}\left(\max_{j \in A} |T_j| > c_n, D \subset \widehat{S}_n\right) \leq \mathbb{P}\left(\max_{j \in A} |\widehat{\beta}_{\widehat{S}_n,j}| > \frac{\widehat{\sigma} c_n}{\sqrt{nC}}, D \subset \widehat{S}_n\right).$$

Let $\gamma = n^{-1} X^T \varepsilon$. Then,

$$\|\widehat{\beta}_A\|^2 \leq \gamma_{\widehat{S}_n}^T Q_{\widehat{S}_n}^2 \gamma_{\widehat{S}_n} \leq \frac{\|\gamma_{\widehat{S}_n}\|^2}{\kappa^2} \leq \frac{k_n \max_{1 \leq j \leq p_n} \gamma_j^2}{\kappa^2}.$$

It follows that
$$\max_{j \in A} |\widehat{\beta}_{\widehat{S}_n, j}| \leq \frac{\sqrt{k_n} \max_{1 \leq j \leq p_n} |\gamma_j|}{\kappa} \leq \sqrt{k_n \log \log n} \max_{1 \leq j \leq p_n} |\gamma_j|,$$
since $\kappa > 0$. So,
$$\mathbb{P}\left(\max_{j \in A} |\widehat{\beta}_{\widehat{S}_n, j}| > \frac{\widehat{\sigma} c_n}{\sqrt{n \log \log n}}, D \subset \widehat{S}_n\right) \leq \mathbb{P}\left(\max_{1 \leq j \leq p_n} |\gamma_j| > \frac{\widehat{\sigma} c_n}{\log \log n \sqrt{n k_n}}\right).$$

Note that $\gamma_j \sim N(0, \sigma^2/n)$, and hence,
$$\mathbb{E}\left(\max_j |\gamma_j|\right) \leq \sqrt{\frac{2\sigma^2 \log(2p_n)}{n}}.$$

There exists $\varepsilon_n \to 0$ such that $\mathbb{P}(B_n) \to 1$ where $B_n = \{(1 - \varepsilon_n) \leq \widehat{\sigma}/\sigma \leq (1 + \varepsilon)\}$. So,
$$\mathbb{P}\left(\max_{1 \leq j \leq p_n} |\gamma_j| > \frac{\widehat{\sigma} c_n}{\log \log n \sqrt{n k_n}}\right) \leq \mathbb{P}\left(\max_{1 \leq j \leq p_n} |\gamma_j| > \frac{\sigma c_n (1 - \varepsilon_n)}{\log \log n \sqrt{n k_n}}, B_n\right)$$
$$\leq \frac{\sqrt{n k_n}}{\sigma (1 - \varepsilon_n) c_n \sqrt{\log \log n}} \mathbb{E}\left(\max_j |\gamma_j|\right)$$
$$\leq \alpha + o(1). \qquad \Box$$

**8. Discussion.** The multi-stage method presented in this paper successfully controls type I error while giving reasonable power. The lasso and stepwise have similar performance. Although theoretical results assume independent data for each of the three stages, simulations suggest that leave-one-out cross-validation leads to valid type I error rates and greater power. Screening the data in one phase of the experiment and cleaning in a followup phase leads to an efficient experimental design. Certainly this approach deserves further theoretical investigation. In particular, the question of optimality is an open question.

The literature on high-dimensional variable selection is growing quickly. The most important deficiency in much of this work, including this paper, is the assumption that the model $Y = X^T \beta + \varepsilon$ is correct. In reality, the model is at best an approximation. It is possible to study linear procedures when the linear model is not assumed to hold as in Greenshtein and Ritov (2004). We discuss this point in the Appendix. Nevertheless, it seems useful to study the problem under the assumption of linearity to gain insight into these methods. Future work should be directed at exploring the robustness of the results when the model is wrong.

Other possible extensions include: dropping the normality of the errors, permitting nonconstant variance, investigating the optimal sample sizes for each stage and considering other screening methods besides cross-validation.

Finally, let us note that the example involving unfaithfulness, that is, cancellations of parameters to make the marginal correlation much different than the regression coefficient, pose a challenge for all the methods and deserve more attention even in cases of small $p$.

## APPENDIX: PREDICTION

Realistically, there is little reason to believe that the linear model is correct. Even if we drop the assumption that the linear model is correct, sparse methods like the lasso can still have good properties as shown in Greenshtein and Ritov (2004). In particular, they showed that the lasso satisfies a risk consistency property. In this appendix we show that this property continues to hold if $\lambda$ is chosen by cross-validation.

The lasso estimator is the minimizer of $\sum_{i=1}^{n}(Y_i - X_i^T \beta)^2 + \lambda \|\beta\|_1$. This is equivalent to minimizing $\sum_{i=1}^{n}(Y_i - X_i^T \beta)^2$ subject to $\|\beta\|_1 \leq \Omega$, for some $\Omega$. (More precisely, the set of estimators as $\lambda$ varies is the same as the set of estimators as $\Omega$ varies.) We use this second version throughout this section.

The predictive risk of a linear predictor $\ell(x) = x^T \beta$ is $R(\beta) = \mathbb{E}(Y - \ell(x))^2$ where $(X, Y)$ denotes a new observation. Let $\gamma = \gamma(\beta) = (-1, \beta_1, \ldots, \beta_p)^T$ and let $\Gamma = \mathbb{E}(ZZ^T)$ where $Z = (Y, X_1, \ldots, X_p)$. Then, we can write $R(\beta) = \gamma^T \Gamma \gamma$. The lasso estimator can now be written as $\widehat{\beta}(\Omega_n) = \arg\min_{\beta \in B(\Omega_n)} \widehat{R}(\beta)$ where $\widehat{R}(\beta) = \gamma^T \widehat{\Gamma} \gamma$ and $\widehat{\Gamma} = n^{-1} \sum_{i=1}^{n} Z_i Z_i^T$.

Define
$$\beta_* = \arg\min_{\beta \in B(\Omega_n)} R(\beta),$$
where
$$B(\Omega_n) = \{\beta : \|\beta\|_1 \leq \Omega_n\}.$$

Thus, $\ell_*(x) = x^T \beta_*$ is the best linear predictor in the set $B(\Omega_n)$. The best linear predictor is well defined even though $\mathbb{E}(Y|X)$ is no longer assumed to be linear. Greenshtein and Ritov (2004) call an estimator $\widehat{\beta}_n$ persistent, or predictive risk consistent, if
$$R(\widehat{\beta}_n) - R(\beta_*) \xrightarrow{P} 0$$
as $n \to \infty$.

The assumptions we make in this section are:

(B1) $p_n \leq e^{n^\xi}$, for some $0 \leq \xi < 1$;

(B2) the elements of $\widehat{\Gamma}$ satisfy an exponential inequality
$$\mathbb{P}(|\widehat{\Gamma}_{jk} - \Gamma_{jk}| > \varepsilon) \leq c_3 e^{-nc_4 \varepsilon^2}$$
for some $c_3, c_4 > 0$;

(B3) there exists $B_0 < \infty$ such that, for all $n$, $\max_{j,k} \mathbb{E}(|Z_j Z_k|) \leq B_0$.

Condition (A2) can easily be deduced from more primitive assumptions as in Greenshtein and Ritov (2004), but for simplicity we take (A2) as an assumption. Let us review one of the results in Greenshtein and Ritov (2004). For the moment, replace (A1) with the assumption that $p_n \leq n^b$, for some $b$. Under these conditions, it follows that

$$\Delta_n \equiv \max_{j,k} |\widehat{\Gamma}_{jk} - \Gamma_{jk}| = O_P\bigg(\sqrt{\frac{\log n}{n}}\bigg).$$

Hence,

$$\sup_{\beta \in B(\Omega_n)} |R(\beta) - \widehat{R}(\beta)| = \sup_{\beta \in B(\Omega_n)} |\gamma^T (\Gamma - \widehat{\Gamma}) \gamma|$$

$$\leq \Delta_n \sup_{\beta \in B(\Omega_n)} \|\gamma\|_1^2 = \Omega_n^2 O_P\bigg(\sqrt{\frac{\log n}{n}}\bigg).$$

The latter term is $o_P(1)$ as long as $\Omega_n = o((n/\log n)^{1/4})$. Thus, we have the following.

THEOREM A.1 [Greenshtein and Ritov (2004)].  *If $\Omega_n = o((n/\log n)^{1/4})$, then the lasso estimator is persistent.*

For future reference, let us state a slightly different version of their result that we will need. We omit the proof.

THEOREM A.2.  *Let $\gamma > 0$ be such that $\xi + \gamma < 1$. Let $\Omega_n = O(n^{(1-\xi-\gamma)/4})$. Then, under* (B1) *and* (B2),

(20) $$\mathbb{P}\bigg(\sup_{\beta \in B(\Omega_n)} |\widehat{R}(\beta) - R(\beta)| > \frac{1}{n^{\gamma/4}}\bigg) = O(e^{-cn^{\gamma/2}})$$

*for some $c > 0$.*

The estimator $\widehat{\beta}(\Omega_n)$ lies on the boundary of the ball $B(\Omega_n)$ and is very sensitive to the exact choice of $\Omega_n$. A potential improvement—and something that reflects actual practice—is to compute the set of lasso estimators $\widehat{\beta}(\ell)$, for $0 \leq \ell \leq \Omega_n$ and then select from that set based on cross validation. We now confirm that the resulting estimator preserves persistence. As before, we split the data into $\mathcal{D}_1$ and $\mathcal{D}_2$. Construct the lasso estimators $\{\widehat{\beta}(\ell) : 0 \leq \ell \leq \Omega_n\}$. Choose $\widehat{\ell}$ by cross validation using $\mathcal{D}_2$. Let $\widehat{\beta} = \widehat{\beta}(\widehat{\ell})$.

THEOREM A.3. *Let $\gamma > 0$ be such that $\xi + \gamma < 1$. Under* (A1), (A2) *and* (A3), *if $\Omega_n = O(n^{(1-\xi-\gamma)/4})$, then the cross-validated lasso estimator $\widehat{\beta}$ is persistent. Moreover,*

(21) $$R(\widehat{\beta}) - \inf_{0 \leq \ell \leq \Omega_n} R(\widehat{\beta}(\ell)) \xrightarrow{P} 0.$$

PROOF. Let $\beta_*(\ell) = \arg\min_{\beta \in B(\ell)} R(\beta)$. Define $h(\ell) = R(\beta_*(\ell))$, $g(\ell) = R(\widehat{\beta}(\ell))$ and $c(\ell) = \widehat{L}(\widehat{\beta}(\ell))$. Note that, for any vector $b$, we can write $R(b) = \tau^2 + b^T \Sigma b - 2b^T \rho$ where $\rho = (\mathbb{E}(YX_1), \ldots, \mathbb{E}(YX_p))^T$.

Clearly, $h$ is monotone nonincreasing on $[0, \Omega_n]$. We claim that $|h(\ell + \delta) - h(\ell)| \leq c\Omega_n \delta$ where $c$ depends only on $\Gamma$. To see this, let $u = \beta_*(\ell)$, $v = \beta_*(\ell + \delta)$ and $a = \ell \beta_*(\ell + \delta)/(\ell + \delta)$ so that $a \in B(\ell)$. Then,

$$\begin{aligned} h(\ell + \delta) &\leq h(\ell) \\ &= R(u) \leq R(a) \\ &= R(v) + R(a) - R(v) \\ &= h(\ell + \delta) + \frac{2\delta}{\ell + \delta}\rho^T v - \frac{\delta(2\ell + \delta)}{(\ell + \delta)^2} v^T \Sigma v \\ &\leq h(\ell + \delta) + 2\delta C + \delta(2\Omega_n + \delta)C, \end{aligned}$$

where $C = \max_{j,k} |\Gamma_{j,k}| = O(1)$.

Next, we claim that $g(\ell)$ is Lipschitz on $[0, \Omega_n]$ with probability tending to 1. Let $\widehat{\beta}(\ell) = \arg\min_{\beta \in B(\ell)} \widehat{R}(\beta)$ denote the lasso estimator and set $\widehat{u} = \widehat{\beta}(\ell)$ and $\widehat{v} = \widehat{\beta}(\ell + \delta)$. Let $\varepsilon_n = n^{-\gamma/4}$. From (20), the following chain of equations hold except on a set of exponentially small probability:

$$\begin{aligned} g(\ell + \delta) = R(\widehat{v}) &\leq \widehat{R}(\widehat{v}) + \varepsilon_n \leq \widehat{R}(v) + \varepsilon_n \\ &\leq R(v) + 2\varepsilon_n = h(\ell + \delta) + 2\varepsilon_n \\ &\leq h(\ell) + c\Omega_n \delta + 2\varepsilon_n = R(u) + c\Omega_n \delta + 2\varepsilon_n \\ &\leq R(\widehat{u}) + c\Omega_n \delta + 2\varepsilon_n = g(\ell) + c\Omega_n \delta + 2\varepsilon_n. \end{aligned}$$

A similar argument can be applied in the other direction. Conclude that

$$|g(\ell + \delta) - g(\ell)| \leq c\Omega_n \delta + 2\varepsilon_n$$

except on a set of small probability.

Now let $A = \{0, \delta, 2\delta, \ldots, m\delta\}$ where $m$ is the smallest integer such that $m\delta \geq \Omega_n$. Thus, $m \sim \Omega_n/\delta_n$. Choose $\delta = \delta_n = n^{-3(1-\xi-\gamma)/8}$. Then, $\Omega_n \delta_n \to 0$ and $\Omega_n/\delta_n \leq n^{3(1-\xi-\gamma)/4}$. Using the same argument as in the proof of Theorem 3.2,

$$\max_{\ell \in A} |\widehat{L}(\widehat{\beta}(\ell)) - R(\widehat{\beta}(\ell))| = \sigma_n,$$

where $\sigma_n = o_P(1)$. Then,

$$\begin{aligned}
R(\beta_*(\Omega_n)) \leq R(\widehat{\beta}) &\leq \widehat{L}(\widehat{\beta}(\widehat{\ell})) + \sigma_n \\
&\leq \widehat{L}(m\delta_n) + \sigma_n \leq g(m\delta_n) + 2\sigma_n \leq g(\Omega_n) + 2\sigma_n + c\Omega_n\delta_n \\
&\leq h(\Omega_n) + 2\sigma_n + \varepsilon_n + c\Omega_n\delta_n \\
&= R(\beta_*(\Omega_n)) + 2\sigma_n + \varepsilon_n + c\Omega_n\delta_n
\end{aligned}$$

and persistence follows. To show the second result, let $\widetilde{\beta} = \arg\min_{0 \leq \ell \leq \Omega_n} g(\ell)$ and $\overline{\beta} = \arg\min_{\ell \in A} g(\ell)$. Then,

$$\begin{aligned}
R(\widetilde{\beta}) &\leq \widehat{L}(\widetilde{\beta}) + \sigma_n \leq \widehat{L}(\overline{\beta}) + \sigma_n \\
&\leq R(\overline{\beta}) + 2\sigma_n \leq R(\widetilde{\beta}) + 2\sigma_n + c\delta_n\Omega_n
\end{aligned}$$

and the claim follows. $\square$

**Acknowledgments.** The authors are grateful for the use of a portion of the sample from the Osteoporotic Fractures in Men (MrOS) Study to illustrate their methodology. MrOs is supported by the National Institute of Arthritis and Musculoskeletal and Skin Diseases (NIAMS), the National Institute on Aging (NIA) and the National Cancer Institute (NCI) through Grants U01 AR45580, U01 AR45614, U01 AR45632, U01 AR45647, U01 AR45654, U01 AR45583, U01 AG18197 and M01 RR000334. Genetic analyses in MrOS were supported by R01-AR051124. We also thank two referees and an AE for helpful suggestions.

Department of Statistics  
Carnegie Mellon University  
Pittsburgh  
USA  
E-mail: larry@stat.cmu.edu  
  roeder@stat.cmu.edu